

\def\N{{\cal N}}
\baselineskip=14pt
\parskip=10pt
\def\halmos{\hbox{\vrule height0.15cm width0.01cm\vbox{\hrule height
  0.01cm width0.2cm \vskip0.15cm \hrule height 0.01cm width0.2cm}\vrule
  height0.15cm width 0.01cm}}
\font\eightrm=cmr8 
\font\eighttt=cmtt8
\magnification=\magstephalf

\def\1{{\overline{1}}}
\def\2{{\overline{2}}}
\parindent=0pt
\overfullrule=0in
\def\Tilde{\char126\relax}
\def\frac#1#2{{#1 \over #2}}
\bf
\centerline
{
The Number of Inversions and the Major Index of Permutations 
}
\centerline
{
are Asymptotically Joint-Independently-Normal (Second Edition)
}
\rm
\bigskip
\centerline
{Written by:}
\centerline
{ 
{\it Andrew BAXTER}{$^1$} and {\it Doron ZEILBERGER}\footnote{$^1$}
{\eightrm  \raggedright
Department of Mathematics, Rutgers University (New Brunswick),
Hill Center-Busch Campus, 110 Frelinghuysen Rd., Piscataway,
NJ 08854-8019, USA.
{\eighttt [baxter,zeilberg]  at math dot rutgers dot edu} ,
\hfill \break
{\eighttt http://www.math.rutgers.edu/\~{}[baxter,zeilberg]} .
This second edition published: Feb. 4, 2011.
First edition was published: April 5, 2010. 
Accompanied by the Maple package {\eighttt InvMaj}
downloadable from
\hfill\break
{\eighttt http://www.math.rutgers.edu/\~{}zeilberg/mamarim/mamarimhtml/invmaj.html} .
\hfill\break
The work of both authors was supported in part by the USA National Science Foundation.
\break
Exclusively published in the Personal Journal of Shalosh B. Ekhad and Doron Zeilberger
\hfill\break
{\eighttt http://www.math.rutgers.edu/\~{}zeilberg/pj.html}  and {\tt arxiv.org}.
DZ is hereby offering \$1000 to the first person to point out a serious
flaw in the argument, that would irreparably invalidate the proof.
All the meta-mathematical opinions expressed here are those of DZ.
}
}

\centerline
{ \eightrm Refereed by:}
\centerline
{\eightrm Mireille Bousquet-M\'elou, Guoniu Han, Emilie Hogan, Svante Janson, Ilias Kotsireas, }
\centerline
{ \eightrm Christian Krattenthaler,  Dan Romik, Vince Vatter, and Herbert Wilf}

{\bf Abstract:} We use {\it recurrences} (alias {\it difference equations}) to prove
that the two most important {\it permutation statistics},
namely the {\it number of inversions} and the {\it major index}, are asymptotically joint-independently-normal.
We even derive  more-precise-than-needed asymptotic formulas for the (normalized) mixed moments.

{\bf Preface to the Second Edition}

In addition to the considerable interest of the results proved in this 
article, and the
even greater interest of the {\it methodology}, this article is a landmark case
in scholarly publishing.
After the first version of this article was {\it outright} 
rejected by an
{\it anonymous} referee 
of the {\it Proceedings of the American Mathematical Society},
because {\it too many details were left to the reader} 
(he or she didn't give us a chance to
write a new version with more details), we decided to 
solicit {\it nine} non-anonymous
reports from world-class experts, assigning them specific parts. 
The division of
labor (quite a few of them did over and above of what we asked them to, 
and refereed everything),
and the full reports, on the first edition, can be gotten from the webpage of this article:

{\eighttt http://www.math.rutgers.edu/\~{}zeilberg/mamarim/mamarimhtml/invmaj.html} ,

already mentioned in footnote 1. With a few exceptions of stylistic suggestions that we preferred
not to adopt, we incorporated {\it all} of their (excellent!) suggestions.
The first edition is also availabe there, for the record.

In order to be faithful to the original version, we have included in square brackets, and
{\eightrm smaller font}, the extra explanations demanded by the referees. We thank them profusely, and
we now believe that the {\it formal correctness} and {\it clarity} far exceeds \%99.99 
of the articles published in (anonymous!) ``peer''-reviewed mathematical journals.
Such journals, even electronic ones,
will soon become obsolete, together with their 
{\it pompous} ``editors'' and
{\it anonymous} referees. 
Instead, the present model of author(s)-appointed
refereeing and self-publishing, in 
the authors' personal websites and the {\it arxiv},
with all the referee reports made public, and the referees
acknowledged and given explicit recognition for their trouble,
would become the norm, possibly with some tweaking.
The present article is {\it exclusively}
published in the Personal Journal of Shalosh B. Ekhad and Doron Zeilberger
and {\tt arxiv.org}.

{\bf Human Statistics}

{\it Human} statistics are numerical attributes defined on {\it humans}, for example,
{\it longevity}, {\it height}, {\it weight}, {\it IQ}, and it is well-known, at least empirically, that
these are, each separately, {\it asymptotically normal}, which means that if you draw a histogram
with the statistical data, it would look like a {\it bell-curve}. It is also true that
they are usually joint-asymptotically-normal, but usually not {\it independently} so.
But if you compute empirically the {\it correlation matrix}, you would get, asymptotically
(i.e. for ``large'' populations) that they are close to being distributed according
to a multivariate (generalized) {\it Gaussian} $exp (-Q(x_1, x_2, ...))$ with $Q(x_1,x_2, \dots)$ a
certain {\it quadratic form} that can be deduced from the correlation matrix.

{\bf Permutation Statistics}

Let our population be the set of {\it permutations} of $\{1,2, \dots, n \}$.
They too, can be assigned {\it numerical attributes}, and the great
classical combinatorialist Dominique Foata (who got his {\it Doctorat de troisi\`eme cycle} in statistics!)
coined the term {\it permutation statistics} for them.

The most important permutation statistic is the {\bf number of inversions}, $inv(\pi)$,
that counts the number of pairs $1 \leq i < j \leq n$ such that $\pi_i>\pi_j$
(and ranges from $0$ to $n(n-1)/2$). For example, $inv(314625)=5$,
corresponding to the set of pairs $\{[1,2],[1,5],[3,5],[4,5],[4,6]\}$.
It features in the definition of the determinant, and Netto
proved that the {\it probability generating function} 
(the polynomial in $q$ such that its coefficient of $q^i$ is the probability that
a uniformly-at-random $n$-permutation has $i$ inversions) is given by
$$
{
{(1)(1+q)(1+q+q^2) \cdots (1+q+q^2+ \dots + q^{n-1} )}
\over
{n!}
}=
{
{\prod_{i=1}^{n} (1-q^i)}
\over
{n!(1-q)^n}
} \quad .
$$

The {\it second most} important permutation statistic is the {\it major index}, $maj(\pi)$, that is the
sum of the places $i$, where $\pi_{i}>\pi_{i+1}$.
For example, $maj(314625)=1+4=5$, because at $i=1$ and $i=4$ we have  descents.
Major Percy Alexander MacMahon [M] famously proved that the probability generating function
for the major index is also given by that very same formula.
In other words the permutation statistics $inv$ and $maj$ are {\it equidistributed}.
Dominique Foata [Fo] gave a lovely seminal {\it bijective} proof that proved the stronger statement
that $inv$ and $maj$ are equi-distributed also when restricted to permutations ending at a given integer.

William Feller ([Fe], $3$rd ed., p.257) proved that the number of inversions (and hence also the major index) is
{\it asymptotically normal} in the following sense. Feller easily computed the {\it expectation},
$$
E[inv]=m_n=n(n-1)/4 \quad ,
$$
and the {\it variance},
$$
\sigma_n^2={{2n^3+3n^2-5n} \over {72}} \quad .
$$
If we denote by $X_n$ the {\it centralized} and {\it normalized} random variable
$$
X_n ={{inv-m_n} \over {\sigma_n}} \quad ,
$$
then $X_n \rightarrow \N$, as $n \rightarrow \infty$, in distribution,
where $\N$ is the Gaussian distribution whose {\it probability density function} is
$e^{-x^2/2}/\sqrt{2\pi}$.

A {\it computer-generated} proof, that gives much more detail about the
rate of convergence to $\N$, can be obtained using Zeilberger's 
Maple package \hfill \break
{\tt http://www.math.rutgers.edu/\Tilde zeilberg/tokhniot/AsymptoticMoments}
that accompanies the article [Z].

So both $inv$ and $maj$ are {\it individually} asymptotically normal, but what about
their {\it interaction}? 
In this article, we prove that they are asymptotically {\it joint-independently-normal}.
In other words, defining,
$$
X_n(\pi):=\frac{inv(\pi)-m_n}{\sigma_n} \quad , \quad  Y_n(\pi):=\frac{maj(\pi)-m_n}{\sigma_n} \quad ,
$$ 
we have that
$$
Pr(X_n \leq s \,\, , \,\, Y_n \leq t )
\rightarrow
\frac{1}{2\pi} \int_{-\infty}^s \int_{-\infty}^t e^{-x^2/2-y^2/2} \, dy \, dx
\quad as \quad n \rightarrow \infty .
$$

{\bf A Brief History}

It all started when the great Swedish probabilist Svante Janson (member of the Swedish Academy of Science,
that awards the Nobel prizes) asked Donald Knuth (one of the greatest computer scientists of all time,
winner of the Turing and Kyoto prizes, among many other honors) about
the {\it asymptotic covariance} of $inv$ and $maj$. Neither of these luminaries knew the answer,
so Don Knuth asked one of us (DZ). DZ didn't know the answer either, so he asked his
beloved servant, Shalosh B. Ekhad, who {\it immediately} ([E]) produced, not just the asymptotics,
but the {\it exact answer}! It turned out to be $n(n-1)/8$. In particular, the
{\it correlation coefficient}, $Cov(inv,maj)/\sigma_n^2= \frac{\frac{n(n-1)}{8}}{\frac{2n^3+3n^2-5n}{72}}=\frac{9}{2n}+O(1/n^2)$
tends to zero as $n$ goes to infinity. It followed that in the long-run, $inv$ and $maj$ are practically {\it uncorrelated}.

But there are lots of pairs of random variables that are uncorrelated yet not independent. 
A convenient way to 
prove that $X_n$ and $Y_n$ are asymptotically independent (we already know that they
are both normal) is to use the {method of moments}, and to prove that the {\it mixed moments}
$$
M_{r,s}(n):=E[X_n^{r} Y_n^{s}] \quad,
$$
tend to the mixed moments of $\N \times \N$, as $n \rightarrow \infty$ . In other words, for $r,s \geq 1$:
$$
\lim_{n \rightarrow \infty} M_{2r,2s}(n)={{(2r)!} \over {2^r r!}}{{(2s)!} \over {2^s s!}} \quad ,
\eqno(EE)
$$
$$
\lim_{n \rightarrow \infty} M_{2r-1,2s}(n) \, = \, 0 \quad ,
\eqno(OE)
$$
$$
\lim_{n \rightarrow \infty} M_{2r,2s-1}(n) \, = \, 0 \quad ,
\eqno(EO)
$$
$$
\lim_{n \rightarrow \infty} M_{2r-1,2s-1}(n) \, = \, 0 \quad .
\eqno(OO)
$$

Ekhad's brilliant approach to the Janson-Knuth question merely settled the case $r=1,s=1$ of
$(OO)$. Of course, because of symmetry $(OE)$ and $(EO)$ are trivially true (before taking the limits!,
i.e. $M_{2r,2s-1}(n) \equiv 0 $ and $M_{2r-1,2s}(n) \equiv 0$).

[{\eightrm
Indeed
if $com([\pi_1, \dots, \pi_n]):=[n+1-\pi_1, \dots, n+1-\pi_n]$ then trivially $maj(\pi)+maj(com(\pi))=n(n-1)/2$ and
$inv(\pi)+inv(com(\pi))=n(n-1)/2$. So $Pr(\{X_n=i,Y_n=j\})=Pr(\{X_n=-i,Y_n=-j\})$
for all $i,j$, and it follows that 
$$
M_{2r,2s-1}(n)=\sum_{i,j} i^{2r}j^{2s-1} Pr(\{X_n=i,Y_n=j\})=
\sum_{i,j} i^{2r}j^{2s-1} Pr(\{X_n=-i,Y_n=-j\})=
$$
$$
\sum_{i,j} (-i)^{2r}(-j)^{2s-1} Pr(\{X_n=i,Y_n=j\})=
-\sum_{i,j} (i)^{2r}(j)^{2s-1} Pr(\{X_n=i,Y_n=j\})=-M_{2r,2s-1}(n) \quad,
$$
so $M_{2r,2s-1}(n)$ equals to its negative, and so must vanish. The proof that
$M_{2r-1,2s}(n)=0$ is similar.}]

One natural approach would be to extend Ekhad's brilliant derivation of $M_{1,1}(n)$ to the general case,
and try to derive closed-form expressions for $M_{r,s}(n)$ for larger $r$ and $s$.
Since Ekhad's proof [E] is so brief, we can cite it here in full.

``Svante Janson asked Don Knuth,
who asked me, about the covariance of $inv$ and $maj$.
The answer is ${{n} \choose {2}}/4$. To prove it, I asked
Shalosh to compute the average of the quantity 
$(inv ( \pi) - E(inv))(maj( \pi) - E(maj)) $ over all permutations of
a given length $n$, and it gave me,
for $n=1,2,3,4,5$, the values $0,1/4,3/4,3/2,5/2$, respectively. Since
we know a priori\footnote{$^2$}
{\eightrm This is the old trick to compute moments
of combinatorial `statistics', described nicely in [GKP], section 8.2,
by changing the order of summation. It applies equally well
to covariance. 
Rather than actually carrying out the
gory details, we observe that this is always a polynomial whose
degree is trivial to bound.
[Added in 2nd edition: the referees didn't find this obvious and asked for 
an explanation.
See the bottom of page 8 and the top of page 9 in the present article.].}
that this is a polynomial of degree $\leq 4$, this must be it! \halmos''.

Obviously this brute-brute-force approach would be hopeless
for deriving polynomial expressions for the moments $M_{r,s}(n)$ for
larger $r$ and $s$. As we will soon see, the degree of
the polynomial $M_{2r,2s}(n)$ is $3(r+s)$, so for example,
in order to (rigorously) guess  $M_{10,10}(n)$,
we would need $31$ data points, and we would have to ask our computers to
examine more than $31! > 0.822 \cdot 10^{34}$ permutations.

However, an {\it  inspired}, still ``empirical'' (yet {\bf fully rigorous})
``brute-force'' approach {\it does} work.
The first step would be to have a more efficient way to
compute the moments $M_{r,s}(n)$, for {\it specific} $n$ and {\it specific}
$r$ and $s$. We will do it by first designing an efficient way
to generate the {\it probability generating function}, let's call it $G(n)(p,q)$, for
the pair of statistics $(inv,maj)$. There are beautiful
closed-form expressions for $G(n)(p,1)$ and $G(n)(1,q)$
(the same one, actually, due to Netto and MacMahon, given in page 2), 
but no such closed-form expression seems to exist for
the bi-variate generating function, so the best that we
can hope for is to find a {\it recurrence scheme}.

{\bf A Combinatorial Interlude}

Let us forget about probability for a few moments, and focus on a fast
algorithm for computing
$$
H(n)(p,q):=\sum_{\pi \in S_n} p^{inv(\pi)}q^{maj(\pi)} \quad ,
$$
for $n$ up to, say, $n=50$.

[
{\eightrm
Referee Dan Romik
believe that we should mention, at this point,
the ``explicit'' formula of Roselle[R] (mentioned by Knuth[K])
in terms of a certain infinite double product for the $q$-exponential
generating function of $H(n)(p,q)$. Romik believes that
this may lead to an alternative proof, that would even imply
a stronger result (a local limit law). 
We strongly doubt this, and
DZ is hereby offering \$300 for the first person to supply such a proof, 
whose length should not exceed the length of this article.

Referee Christian Krattenthaler believes that we should also mention
the beautiful extension of Roselle's result, 
by Adriano Garsia and Ira Gessel [GG],
handling more permutation statistics.
}
]

Define the {\it weight} of a permutation 
$\pi$ to be $p^{inv(\pi)}q^{maj(\pi)}$.
Suppose that $\pi \in S_n$ ends with $i$, so we can write
$\pi=\pi' i$, where $\pi'$ is a permutation of 
$\{1, \dots , i-1, i+1, \dots n \}$.

When you chop off $i$ from $\pi$ to form $\pi'$ you always lose $n-i$
inversions (that is, $\pi'$ has $n-i$ fewer inversions than $\pi$).  The
major index, however, decreases by $n-1$ if the last letter of $\pi'$, let's
call it $j$, is larger than $i$. If $j<i$ the major index does not change at
all.  So writing $\pi=\pi'' j i$, we have
$$
inv(\pi'' ji)=inv(\pi''j)+n-i \quad,
$$
$$
maj(\pi'' ji) = \cases{
maj(\pi''j) ,& if $j<i$ ;\cr
maj(\pi''j)+n-1 ,&  if $j>i$.\cr}   
$$
Combining, we have
$$
weight(\pi'' ji) = \cases{
p^{n-i}weight(\pi''j) ,& if $j<i$ ;\cr
p^{n-i}q^{n-1} weight(\pi''j),&  if $j>i$.\cr}  
$$
So in order to compute $H(n)(p,q)$, we need to introduce the
more general weight-enumerators of those permutations in $S_n$
that end with an $i$. Let's call these $F(n,i)(p,q)$.
In symbols:
$$
F(n,i)(p,q):=\sum_
{
{{\pi \in S_n} \atop {\pi_n=i}}
} p^{inv(\pi)}q^{maj(\pi)} \quad .
$$
It follows that (let's omit the arguments $(p,q)$ from now on):
$$
F(n,i)=
p^{n-i} \sum_{j=1}^{i-1} F(n-1,j)+p^{n-i}q^{n-1} \sum_{j=i}^{n-1} F(n-1,j) 
\quad .
\eqno(Fni)
$$
  Note that, when we chop off the last entry, $i$, from $\pi=\pi'' j i$, 
$\pi'' j$ is a permutation of $\{1, \ldots, i-1, i+1, \ldots, n\}$.  We
then ``reduce'' $\pi''j$ to a permutation of $\{1,\ldots, n-1\}$ by
diminishing all entries larger than $i$ by 1.  Hence the summation ranges from
$j=i$ to $j=n-1$ rather than from $j=i+1$ to $j=n$.

Replacing $i$ by $i+1$ in the above equation, we have:
$$
F(n,i+1)=
p^{n-i-1} \sum_{j=1}^{i} F(n-1,j)+p^{n-i-1}q^{n-1} 
\sum_{j=i+1}^{n-1} F(n-1,j) \quad .
$$
Subtracting the former equation from $p$ times the latter we get
$$
F(n,i)-pF(n,i+1)=
$$
$$
p^{n-i} \sum_{j=1}^{i-1} F(n-1,j)+p^{n-i}q^{n-1} \sum_{j=i}^{n-1} F(n-1,j) 
$$
$$
-p^{n-i} \sum_{j=1}^{i} F(n-1,j)-p^{n-i}q^{n-1} \sum_{j=i+1}^{n-1} F(n-1,j)
$$
$$
=p^{n-i} \left ( \, \sum_{j=1}^{i-1} F(n-1,j) - \sum_{j=1}^{i} F(n-1,j) \, \right )
+p^{n-i}q^{n-1} \left ( \, \sum_{j=i}^{n-1} F(n-1,j)  - \sum_{j=i+1}^{n-1} F(n-1,j) \,  \right )
$$
$$
=-p^{n-i}F(n-1,i)+p^{n-i}q^{n-1}F(n-1,i)=p^{n-i}(q^{n-1}-1)F(n-1,i) \quad.
$$
Rearranging, we get:
$$
F(n,i)=pF(n,i+1)+p^{n-i}(q^{n-1}-1)F(n-1,i) \quad for \quad  1 \leq i < n \quad .
\eqno(RecF)
$$
We still need to specify $F(n,n)$, and for this
we do need the $\sum$ symbol, namely we use Eq. $(Fni)$ with
$i=n$:
$$
F(n,n)= \sum_{j=1}^{n-1} F(n-1,j)
\quad .
\eqno(Fnn)
$$
The recurrence $(RecF)$ together with the final condition $(Fnn)$,
and the trivial {\it initial condition} $F(1,i)=\delta_{i,1}$,
enables us to efficiently compute $F(n,i)$ for  numeric $(n,i)$, 
for $\{(n,i) \vert 1 \leq i \leq n \leq N\}$ 
for any finite $N$ (not too large, but not too small either: e.g., $N=100$ is still feasible).
In particular, we can compile a table of $H(n)(p,q)=F(n+1,n+1)(p,q)$, 
(the generating function for {\it all} $n$-permutations) for $n \leq N-1$.

\vfill
\eject
{\bf A crash course in multivariable enumerative probability}

Suppose that you have a finite set of objects $S$ and
several {\it statistics} $f_1(s), \dots , f_r(s)$.
The {\it multi-variable generating function}
(weight-enumerator under the weight $x_1^{f_1(s)} \cdots x_r^{f_r(s)}$)
is defined to be:
$$
\sum_{s \in S} x_1^{f_1(s)} \cdots x_r^{f_r(s)} \quad .
$$
Suppose that you pick an element $s \in S$ {\it uniformly at random}
and you want the multivariable generating function such that
the coefficient of $x_1^{a_1} \cdots x_r^{a_r}$ would give you the
probability that $f_1(s)=a_1, \dots, f_r(s)=a_r$.
It is given by:
$$
P(x_1, \dots, x_r)=
{{1} \over {|S|}}\sum_{s \in S} x_1^{f_1(s)} \cdots x_r^{f_r(s)} \quad .
$$
The {\it expectations}, $\bar f_1, \dots, \bar f_r$ are simply
$$
\bar f_i =\left ( {{\partial} \over {\partial x_i}} P \right) (1, \dots, 1) 
\quad.
$$
The {\it centralized} probability generating function is
$$
\tilde P(x_1, \dots, x_r)=\frac{P(x_1, \dots, x_r)} 
{x_1^{\bar f_1}  \cdots x_r^{\bar f_r}} \quad .
$$
The mixed moments (about the mean) of the statistics
$f_1(s), \dots, f_r(s)$ are defined by
$$
Mom[a_1, \dots, a_r]=
{{1} \over {|S|}}
\sum_{s \in S} (f_1(s)-\bar f_1)^{a_1} \cdots (f_r(s)-\bar f_r)^{a_r} \quad .
$$
Often it is more convenient to consider the
mixed {\it factorial moments}, using the combinatorial ``powers''
$z^{(r)}:=z(z-1) \cdots (z-r+1)$, better known as the
{\it falling-factorials}. The mixed factorial moments
are defined analogously by
$$
FM(a_1, \dots, a_r)=
{{1} \over {|S|}}
\sum_{s \in S} (f_1(s)-\bar f_1)^{(a_1)} \cdots (f_r(s)-\bar f_r)^{(a_r)} 
\quad .
$$
Once you know the $FM$'s for all $a_1, \dots, a_r \leq M$, you can
easily figure out the $Mom$'s (for $a_1, \dots, a_r \leq M$), 
using Stirling numbers of the
second kind (see [Z]). 
It is well-known and easy to see that one can just as well use the
method of factorial moments in order to prove asymptotic independence.
In other words, it would suffice to prove the analogs of $(EE),(OE),(EO),(OO)$
with the moments $E[X^rY^s]$ replaced by the factorial moments 
$E[X^{(r)}Y^{(s)}]$.

The best way to compute $FM(a_1, \dots, a_r)$
for all $0 \leq a_1, \dots, a_r \leq M$, for some fixed positive integer
$M$, is via the {\it Taylor expansion}  of $\tilde P(x_1, \dots, x_r)$
around $(x_1, \dots, x_r)=(1,\dots, 1)$, or equivalently,
the {\it Maclaurin expansion} of 
$\tilde P(1+x_1, \dots, 1+x_r)$
$$
\tilde P(1+x_1, \dots, 1+x_r)=\sum_{\alpha_1, \dots, \alpha_r \geq 0}
\frac{FM(\alpha_1, \dots, \alpha_r)}{\alpha_1 ! \cdots \alpha_r !}
 x_1^{\alpha_1} \cdots x_r^{\alpha_r} \quad .
$$

{\bf Back to inv-maj}

In the present approach, we need to put up with the
more general discrete function $F(n,i)(p,q)$ even though
ultimately we are only interested in $H(n)(p,q)=F(n+1,n+1)(p,q)$.
We will prove the stronger statement that even if you restrict
attention to those $(n-1)!$ permutations that end with a specific
$i$, the pair $(inv,maj)$ is still asymptotically-joint-independently normal.

Since the averages of both $inv$ and $maj$ over the permutations that end in $i$ is $n-i+(n-1)(n-2)/4$
[{\eightrm for $inv$ it is obvious, the last entry $i$ contributes $n-i$ to the number of inversions, and
removing the last entry yields an $n-1$-permutation, and for $maj$ this follows from 
Foata's bijection mentioned above that maps $inv$ to $maj$ preserving the last entry}],
the {\it centralized probability generating function} corresponding
to $F(n,i)(p,q)$ is:
$$
G(n,i)(p,q):=
\frac{F(n,i)(p,q)}{(n-1)!(pq)^{n-i+(n-1)(n-2)/4}} \quad .
$$
The recurrence $(RecF)$ becomes
$$
G(n,i)={{1} \over {q}}G(n,i+1)+
\frac{p^{n-i}(q^{n-1}-1)}{(pq)^{n/2} (n-1)}G(n-1,i) \quad ,
\eqno(RecG)
$$
and the {\it final condition} becomes
$$
G(n,n)= \frac{1}{n-1} \sum_{j=1}^{n-1} (pq)^{n/2-j} G(n-1,j)\quad .
\eqno(Gnn)
$$
We also need the obvious {\it initial condition}
$G(1,i)=\delta_{i,1}$.

{\bf Guessing Polynomial Expressions for the Factorial Moments}

  Equipped with these efficient recurrences, our computer computes
$G(n,i)(p,q)$ for many values of $n$ and $i$.  Then for each of these it
computes the $(r,s)$ (mixed) factorial moments $FM(r,s)(n,i)$, for many small numeric
$(r,s)$ by computing the (initial terms of the) Maclaurin series for
$G(n,i)(1+p,1+q)$.  We then fix numeric values of $(r,s)$ and use {\it
polynomial interpolation} to {\it guess} explicit polynomial expressions for
$FM(r,s)(n,i)$ as polynomials in $(n,i)$ for that pair $(r,s)$.  The process
is repeated for all pairs $(r,s)$ for which $0 \leq r,s \leq M$,
for some pre-determined specific positive integer $M$.
We note that it is obvious, both from the 
combinatorics and from the recurrences, that the $FM(r,s)(n,i)$ are always 
polynomials in $(n,i)$, for any fixed numeric $r$ and $s$.

[{\eightrm As we have already mentioned in footnote 2, 
most of the referees didn't find this obvious. 
The proof via the recurrences $(RecG')$ and $(Gnn')$ to be
derived in page 10 is by induction on $(r,s)$ and the fact that
the indefinite sum of a polynomial (in this case with respect to $i$)
is yet another polynomial, and the ``operator'' on the
right of $(Gnn')$ is polynomial-preserving.

Let's sketch the combinatorial proof that the mixed moments $M_{r,s}(n)$
are polynomials in $n$. The combinatorial proof that  $FM(r,s)(n,i)$ are polynomials
in both $n$ and $i$ is similar.
Write $inv(\pi)$ and $maj(\pi)$ as a sum of ``atomic'' events, e.g. for
$inv$ the sum of $\chi(\pi_j>\pi_i)$ over all pairs of integers
$(i,j)$ satisfying $1 \leq i < j \leq n$. Here $\chi(S)=1$ if
$S$ is true and $\chi(S)=0$ if it is false. $maj$ can be similarly expressed
as a sum of $\chi(i \leq j \,\, AND \,\, \pi_j > \pi_{j+1})$.
The sum of $inv(\pi)^r maj(\pi)^s$ 
over all permutations $\pi \in S_n$ 
can be expressed as a multi-sum with 
the outer sum ranging over $S_n$ and the inner sums with
$2(r+s)$ sigma signs. 
Now do discrete Fubini! Bring the formerly outer-sigma, over $S_n$,
{\it all the way inside} past all the other $2(r+s)$ sigma signs.
Each individual sigma sign involves two indices, and collectively
these $2(r+s)$ sigma signs involve $\leq 2(r+s)$ indices, corresponding
to locations in a generic  $n$-permutation $\pi$,
some of whom may coincide.
Let's call them $i_1 < ... <i_k$ ( where $k \leq 2r+2s$).
These can be placed in lots of possible intertwining ways, 
and so can the values of $\pi$ in those places.
There are still finitely many scenarios.
(Formally, one gets a Cartesian product of two 
partially ordered sets, one for the domain and one for the range,
each of which
has finitely many linear extensions. This is reminiscent of
Richard Stanley's theory of P-partitions).
For each particular such scenario
(linear extension of the domain-poset and the range-poset), 
there are ${{n} \choose {k}}$
ways to choose the participant indices,
and ${{n} \choose {k}}$ ways to choose their occupants (i.e. the values of $\pi$ there)
and the remaining $n-k$ entries can, of course, 
be arranged in $(n-k)!$ ways,
yielding ${{n} \choose {k}}^2(n-k)!$ ways. Dividing by $n!$  gives
${{n} \choose {k}}^2(n-k)!/n!={{n} \choose {k}}/k!$,
a polynomial in $n$ of degree $k$. 
Now the whole thing is a sum of
finitely many such ${{n} \choose {k}}/k!$ for $0 \leq k \leq 2r+2s$,
and since a finite sum of polynomials is still a polynomial, we are done!
Now isn't that {\it obvious}?!}
]

It would have been nice if we could guess {\it closed-form} expressions
for $FM(r,s)(n,i)$ for {\it symbolic} $(r,s)$, 
but no such closed-form exists as
far as we know, and besides it is too much to ask for and more than we need.
To prove {\it asymptotic normality} we only need the
{\it leading terms}. Viewing the leading terms, our beloved computer
easily conjectures the following expressions.
For integers $r \geq 0, s \geq 0$ we have:
$$
FM(2r,2s)(n,i)=\frac{(2r)!}{2^r r!}\frac{(2s)!}{2^s s!}
\left ( \frac{1}{36} \right )^{r+s} n^{3r+3s}
\, + \, (lower-total-degree-terms-in-(n,i)) \quad .
$$
[{\eightrm Note that the coefficients of $n^{3r+3s-1}i, n^{3r+3s-2}i^2$ etc.
are all zero, hence they don't show up!}].

For integers $r \geq 0, s \geq 1$ we have:
$$
FM(2r,2s-1)(n,i)=\frac{(2r)!}{2^r r!}\frac{(2s)!}{2^s s!}
\left ( \frac{1}{36} \right )^{r+s-1} n^{3r+3s-6}
[-(s-1)n^3-6rn^2i+18rni^2-12ri^3] 
$$
$$	
\, + \, (lower-total-degree-terms-in-(n,i)) \quad .
$$

For integers $r \geq 1, s \geq 0$ we have:
$$
FM(2r-1,2s)(n,i)=-\frac{(2r)!}{2^r r!}\frac{(2s)!}{2^s s!}
\left ( \frac{1}{36} \right )^{r+s-1} (r-1)n^{3r+3s-3}
\, + \, (lower-total-degree-terms-in-(n,i)) \quad .
$$
For integers $r \geq 1, s \geq 1$ we have:
$$
FM(2r-1,2s-1)(n,i)=\frac{(2r)!}{2^r r!}\frac{(2s)!}{2^s s!}
\left ( \frac{1}{36} \right )^{r+s-1} \frac{9}{2} n^{3r+3s-6} (n-2i)^2
\, + \, (lower-total-degree-terms-in-(n,i)) \quad .
$$

{\bf Nice conjectures but what about proofs?}

While we prefer the empirical approach of guessing, an alternative
approach to finding many $FM(r,s)(n,i)$'s,
that is also necessary in order to rigorously prove our conjectures,
is to first use
$(RecG)$ and $(Gnn)$. Write $G(n,i)(1+p,1+q)$ as an {\it infinite}
{\bf generic} Taylor series around $(0,0)$, and write down the
implied {\it infinite-order} recurrences expressing
$FM(r,s)$ in terms of $FM(r',s')$ with $r'+s' <r+s$.

[{\eightrm
The infinite-order recurrence for the $FM(r,s)(n,i)$, obtained from $(RecG)$
is gotten by expanding $1/(1+q)$ and
$$
\frac{(1+p)^{n-i} (1+q)^{n-i}-1}{((1+p)(1+q))^{n/2}} \quad
$$
as Maclaurin series in $(p,q)$, using the binomial theorem and 
manipulations on formal power series (that Maple does automatically to any desired order),
and combining terms.
Similarly, the implication of $(Gnn)$ is obtained by expanding $((1+p)(1+q))^{n/2-j}$ using
the binomial theorem. Both of these tasks are accomplished by procedure
{\tt MOP} in the Maple package {\tt InvMaj}. } ]

Note that in order to compute $FM(r,s)(n,i)$, for any specific, numeric $r$ and $s$,
we only need finitely
many terms (actually $rs-1$ of them) of the infinite-order recurrence, since eventually
all the contributions will be zero.
Of course, as we have already commented, there is
no hope for finding a general expression for 
$FM(2r,2s)(n,i)$ , $FM(2r,2s-1)(n,i)$ , $FM(2r-1,2s)(n,i)$  and
$FM(2r-1,2s-1)(n,i)$, depending {\it explicitly} on $r$ and $s$ (i.e. {\bf symbolically} in terms of $r$ and $s$)
as well as on $n$ and $i$, but to prove, by induction on $r,s$,
that the above {\it leading terms}  are valid, all we need to do is
to verify that the leading terms of the implied recurrences for the
$FM(r,s)$'s (that we have just talked about)
are consistent with the above explicit expressions.

The implication of $(RecG)$ is
$$
FM(r,s)(n,i)-FM(r,s)(n,i+1)=
$$
$$
-sFM(r,s-1)(n,i+1)+sFM(r,s-1)(n-1,i)+ 
\frac{rs(n-2i)}{2}FM(r-1,s-1)(n-1,i) \, + \,
$$
$$
(lower-order-terms) \quad ,
\eqno(RecG')
$$
while the implication of $(Gnn)$ is:
$$
FM(r,s)(n,n)=
\frac{1}{n-1} \sum_{j=1}^{n-1} FM(r,s)(n-1,j) -
\frac{s}{2(n-1)} \sum_{j=1}^{n-1} (2j-n)FM(r,s-1)(n-1,j) -
$$
$$
\frac{r}{2(n-1)} \sum_{j=1}^{n-1} (2j-n)FM(r-1,s)(n-1,j) 
+
\frac{rs}{4(n-1)} \sum_{j=1}^{n-1} (2j-n)^2 FM(r-1,s-1)(n-1,j) 
$$
$$
+ (lower-order-terms) \quad .
\eqno(Gnn')
$$
[{\eightrm
Referee Guoniu Han
correctly commented that one must say something about the
degree of the polynomials in $(n,i)$ that reside in the ``lower-order terms''
that feature as coefficients in these
infinite-order recurrences. It turns out the the coefficient of
$FM(r-r',s-s')(n-1,i)$ in $(RecG')$
and the coefficient of $FM(r-r',s-s')(n,i+1)$ (and similarly for $(Gnn')$) are polynomials
in $(n,i)$ of total degree $\leq r'+s'$, thanks to the 
binomial theorem, and this would imply, by induction on $(r,s)$ that
not only are $FM(2r,2s)(n,i)$ etc. polynomials, but their leading terms look as
claimed above.}]

The next step is to spell out 
these two recurrences each
into the four cases according to whether
$(r,s)$ is $(even, even)$, $(even, odd)$, $(odd, even)$, and $(odd, odd)$.

Once you have these eight recurrences, 
for each and every one of them,
you plug in the above conjectured
expressions for the leading terms, and verify that up to the
leading terms, things agree. At the end of the day, after dividing by 
$$
\frac{(2r-2)!}{(r-1)!2^{r-1}} \frac{(2s-2)!}{(s-1)!2^{s-1}} \quad ,
$$
this boils down 
to proving equalities among certain low-degree polynomials in $(r,s)$ 
(namely the leading coefficients, in $(n,i)$), 
that in turn, reduces (since $A=B$ iff $A-B=0$)
to proving that certain low-degree polynomials in $(r,s)$  
are identically zero.

So in order to check that these low-degree polynomials in $(r,s)$ are all 
identically zero,
it is enough to check each of them for finitely 
(and not-too-many) {\it numeric} $r,s$.
Typing {\tt Check1(FM8m(n,i),n,i);}  and {\tt Check2(FM8m(n,i),n,i);}   
in the Maple package {\tt InvMaj} does exactly that, by checking that if
you plug in the conjectured leading terms of  $FM(2r,2s)(n,i)$,
$FM(2r-1,2s)(n,i)$, $FM(2r,2s-1)(n,i)$ and $FM(2r-1,2s-1)(n,i)$
and subtract the right sides from the left sides 
(for each of the eight cases)
you get  lower-order polynomials, in $(n,i)$, for all $1 \leq r,s \leq 8$.
This proves all these claims (rigorously), with a vengeance! 
The $(8/2)^2=16$ special cases are much more than is needed, since the
relevant polynomials in $(r,s)$ are easily seen to have 
degree $\leq 2$ so $(2+1)^2=9$ agreements would have sufficed.

{\bf The Maple package InvMaj}

All the nitty-gritty calculations described above, that
constitute  {\bf a fully rigorous} proof, may be found
in the Maple package {\tt InvMaj} accompanying this article.
This package is available from the webpage of the present article:

{\tt http://www.math.rutgers.edu/\Tilde zeilberg/mamarim/mamarimhtml/invmaj.html} \quad ,

where the reader can also find some sample input and output.
The direct url of the package is: \hfill \break {\tt http://www.math.rutgers.edu/\Tilde zeilberg/tokhniot/InvMaj} \quad .

{\bf La Grande Finale}

The special cases $r=1,s=0$ and $r=0,s=1$ give
$$
FM(2,0)(n,i)=\frac{1}{36}n^3+ O(n^2) \quad ,
$$
$$
FM(0,2)(n,i)=\frac{1}{36}n^3+ O(n^2) \quad .
$$
So (recall that we are interested in the {\it normalized} mixed factorial moments)
$$
\frac{FM(2r,2s)(n,i)}{FM(2,0)(n,i)^rFM(0,2)(n,i)^s}=
\frac{(2r)!}{2^r r!}\frac{(2s)!}{2^s s!}+O(1/n) \quad ,
$$
$$
\frac{FM(2r,2s-1)(n,i)}{FM(2,0)(n,i)^rFM(0,2)(n,i)^{s-1/2}}=
o(1/n) \quad ,
$$
$$
\frac{FM(2r-1,2s)(n,i)}{FM(2,0)(n,i)^{r-1/2}FM(0,2)(n,i)^s}=
o(1/n) \quad ,
$$
$$
\frac{FM(2r-1,2s-1)(n,i)}{FM(2,0)(n,i)^{r-1/2}FM(0,2)(n,i)^{s-1/2}}=
O(1/n) \quad .
$$
And we see that as $n \rightarrow \infty$ these indeed converge
to the mixed moments of the famous mixed moments
of the bivariate independent normal distribution
$e^{-a^2/2-b^2/2}/(2 \pi)$ \halmos.

{\bf Encore: A more refined asymptotics for the (Normalized) Mixed Moments}

With more effort, we (or rather, our computer) can guess-and-prove
the following asymptotics for the case of interest
$(n+1,n+1)$, i.e. the asymptotic expressions for the
centralized-and-normalized (genuine, not factorial)  mixed-moments,
let's call them $\alpha(r,s)(n)$,
for the pair of random variables $(inv,maj)$ acting on the set of
permutations of length $n$.
Indeed, according to S. B. Ekhad, we have:
$$
\alpha (2r,2s)(n)=\frac{(2r)!}{2^r r!}\frac{(2s)!}{2^s s!}
\left (  1- \frac{9( r^2+s^2-r-s)}{25} 
\cdot \frac{1}{n} +O(\frac{1}{n^2}) \right ) \quad ,
$$
$$
\alpha (2r-1,2s-1)(n)=\frac{(2r)!}{2^r r!}\frac{(2s)!}{2^s s!}
\left ( \frac{9}{2n}+
\left ( -\frac{81}{50}(r^2+s^2)+\frac{243}{50}(r+s)-\frac{1773}{100}
\right )
\frac{1}{n^2}
+O(\frac{1}{n^3}) \right ) \quad .
$$
Of course, by {\it symmetry} $\alpha(2r,2s-1)(n)$ and
$\alpha(2r-1,2s)(n)$ are {\it identically} (not just asymptotically!)
zero. 

{\bf References}

[E] Shalosh B. Ekhad, {\it The joy of brute force: the covariance
of the major index and the number of inversions}, Personal Journal
of S. B. Ekhad and D. Zeilberger, \hfill\break
{\tt http://www.math.rutgers.edu/$\sim$ zeilberg/pj.html}, ca. 1995.

[Fe] William Feller, ``{\it An Introduction to Probability
Theory and Its Application}'', volume 1, three editions.
John Wiley and sons. First edition: 1950. Second edition:
1957. Third edition: 1968.

[Fo] Dominique Foata, {\it On the Netto inversion number of a sequence}, Proc. Amer. Math. Soc. {\bf 19} (1968), 236-240.

[GG] Adriano Garsia and Ira Gessel, {\it  Permutations statistics and partitions}, Adv. in Math. {\bf 31} (1979),
no. 3, 288-305.

[GKP] Ronald Graham,  Donald E. Knuth, and Oren Patashnik,
{\it Concrete Mathematics: A Foundation for
Computer Science}, Addison Wesley, Reading, 1989.

[M] Major P. A. MacMahon, 
{\it The indices of permutations and the derivation
therefrom of functions of a single variable
associated with the permutations of
any assemblage of objects}, Amer. J. Math. {\bf 35} (1913), 281-322.

[K] Donald E. Knuth, {\it ``The Art of Computer Programming, Vol. 3: Sorting and Searching, 2nd Ed.''}, 
Solution to Exercise 27 on page 596.

[R] David P. Roselle, {\it Coefficients associated with the expansion of certain products},
Proc. Amer. Math. Soc. {\bf 45} (1974), 144-150.

[Z] Doron Zeilberger, {\it The automatic central limit theorems
generator (and much more!)}, in: ``Advances in Combinatorial
Mathematics'' (in honor of Georgy P. Egorychev), Ilias S. Kotsireas
and Eugene V. Zima, eds., Springer, 2009.

\end